\pdfoutput=1
\documentclass{shinyart}

\usepackage[utf8]{inputenc}

\usepackage{shinybib}

\usepackage[utf8]{inputenc}

\usepackage{enumitem}
\setlist[enumerate]{label={(\roman*)}}
\usepackage{autonum}

\renewcommand{\phi}{\varphi}
\newcommand{\N}{\mathbb{N}}
\newcommand{\R}{\mathbb{R}}
\newcommand{\dual}[1]{\langle #1 \rangle}

\newcommand{\norm}[1]{\| #1 \|}
\newcommand{\set}[2]{\left\{#1:#2\right\}}
\newcommand{\wkto}{\rightharpoonup}

\DeclareMathOperator{\dom}{\mathrm{dom}}
\DeclareMathOperator{\ran}{\mathrm{ran}}
\DeclareMathOperator{\Id}{\mathrm{Id}}

\newcommand{\cR}{\mathcal{R}}
\newcommand{\ydel}{y^\delta}
\newcommand{\udag}{u^\dagger}
\newcommand{\ukd}{u_{k}^{\delta}}
\newcommand{\ukdn}{u_{k+1}^{\delta}}

\title{Regularization of ill-posed problems with non-negative solutions}
\dedication{\sffamily \large Dedicated to the memory of Jonathan M. Borwein}
\author{%
    Christian Clason\thanks{%
        Faculty of Mathematics, 
        University Duisburg-Essen, 
        45117 Essen, Germany 
    (\email{christian.clason@uni-due.de})}
    \and 
    Barbara Kaltenbacher\thanks{%
        Institute of Mathematics, 
        Alpen-Adria Universität Klagenfurt, 
        Universitätsstrasse 65--67 
        9020 Klagenfurt, Austria 
    (\email{barbara.kaltenbacher@aau.at}, \email{elena.resmerita@aau.at})}
    \and 
    Elena Resmerita\footnotemark[2]
}

\hypersetup{
    pdftitle={Regularization of ill-posed problems with non-negative solutions},
    pdfauthor={Christian Clason, Barbara Kaltenbacher, Elena Resmerita},
}

\begin{document}

\maketitle

\begin{abstract}
    This survey reviews variational and iterative methods for reconstructing non-negative solutions of ill-posed problems in infinite-dimensional spaces.
    We focus on two classes of methods: variational methods based on entropy-minimization or constraints, and iterative methods involving projections or non-negativity-preserving multiplicative updates.
    We summarize known results and point out some open problems.
\end{abstract}

\section{Introduction} 

Many inverse problems are concerned with the reconstruction of parameters that are a priori known to be non-negative, such as material properties or densities (in particular, probability densities). Non-negative solutions also frequently occur in astronomy and optical tomography, in particular in Poisson models for positron emission tomography (PET), see \cite{shepp_vardi82, vardi_shepp_kaufmann85}.
Note that the literature on the finite-dimensional setting is very rich, and quite comprehensive surveys are already available; see, e.g., \cite{byrne2007,Byrne2015a,Byrne2015b}. 

Borwein and collaborators dealt in a series of papers with the case involving  operators with finite-dimensional range; see, e.g., \cite{BorweinLewis91, BorweinLewis91a,BorweinGoodrichLimber96}. Concerned with reconstructing a density function from a finite number $n$ of density moments, they have approached the problem from a few perspectives. For instance, it has been shown in \cite{BorweinLewis91} that the best (Boltzmann--Shannon) entropy estimates converge in the $L^1$-norm to the best entropy estimate of the limiting problem as $n\to\infty$. Along with this result, strong properties of the Boltzmann--Shannon entropy such as a Kadec--Klee property have been derived. Note that the dual problem of the maximum entropy estimates problem has been quite instrumental in showing further results such as error bounds. From a computational point of view, choosing one entropy (e.g., Dirac--Fermi or Burg) over the other has been the main topic in \cite{BorweinGoodrichLimber96}. The work \cite{Borwein93} studies the case of operators with infinite-dimensional range by proposing relaxed problems in the spirit of Morozov and Tikhonov regularization (cf.~\cref{sec:variational}).

The context of infinite-dimensional function spaces for reconstructing non-negative solutions of ill-posed operator equations has been much less investigated in the literature. Therefore, this work focuses on methods for problems in such spaces from a deterministic point of view.

We will primarily consider linear operator equations
\begin{equation}\label{eq:op_equation}
    Au=y
\end{equation}
with the operator $A:X\to Y$ mapping between suitable infinite-dimensional function spaces $X$ and $Y$. We assume that \eqref{eq:op_equation} admits a non-negative solution $u^\dag \geq 0$ (which we will make precise below) and that it is ill-posed in the sense that small perturbations of $y$ can lead to arbitrarily large perturbations on $u$ (or even lead to non-existence of a solution). Besides enforcing non-negativity of the solution for given data, a solution approach therefore also needs to have regularizing properties, i.e., be stable even for noisy data $y^\delta$ with 
\begin{equation}\label{eq:delta}
    \|\ydel-y\|_Y\leq\delta
\end{equation}
in place of $y$ and yield reconstructions $u^\delta$ that converge to $u^\dag$ as $\delta\to 0$. Two approaches are wide-spread in the literature:
\begin{enumerate}
    \item \emph{Variational methods} are based on minimizing a weighted sum of a discrepancy term and a suitable regularization term (\emph{Tikhonov regularization}) or on minimizing one of these terms under a constraint on the other (\emph{Ivanov} or \emph{Morozov regularization}, respectively).
    \item \emph{Iterative methods} construct a sequence of iterates approximating -- for exact data -- the solution $u^\dag$; regularization is introduced by stopping the iteration based on a suitable discrepancy principle.
\end{enumerate}
Regarding the regularization theory for ill-posed problems, we refer, e.g., to the classical work \cite{Engl}; of particular relevance in the context of non-negative solutions are regularization terms or iterations based on the Boltzmann--Shannon entropy and the associated Kullback--Leibler divergence, and we will focus especially on such methods.

\bigskip

The manuscript is organized as follows. \Cref{sec:preliminaries} recalls useful algebraic and topological properties of the mentioned entropy functionals. \Cref{sec:variational} reviews several variational entropy-based regularization methods (Morozov, Tikhonov, Ivanov), while \cref{sec:iterative}
is dedicated to iterative methods for general linear ill-posed equations, both ones involving projections onto the non-negative cone and ones based on multiplicative updates preserving non-negativity.

\section{Preliminaries}
\label{sec:preliminaries}

Let $\Omega$ be an open and bounded subset of ${\mathbb{R}}^d$.
The \emph{negative of the Boltzmann--Shannon entropy} is the function
$f:L^1(\Omega)\rightarrow (-\infty,+\infty]$,
given by\footnote{We use the convention $0\log 0=0$.}
\begin{equation}\label{eq:entropy}
    f(u)=
    \begin{cases}
        \int_{\Omega} u(t)\log u(t)\,dt & \text{if $u\in L^1_+(\Omega)$ and $u\log u\in {{L}}^1(\Omega)$}, \\
        {+\infty} & \text{otherwise.}
    \end{cases}
\end{equation}
Here and in what follows, we set for $p\in[1,\infty]$
\begin{equation}
    L_+^p(\Omega) := \set{u\in L^p(\Omega)}{u(x) \geq 0 \text{ for almost every }x\in\Omega},
\end{equation}
while $\|\cdot\|_p$ denotes, as usual, the norm of the space $L^p(\Omega)$.

We recall some useful properties of the negative Boltzmann--Shannon entropy from, e.g., \cite[proof of Thm.~1]{amato_hughes1991}, \cite[Lem.~2.1, 2.3]{eggermont1993}, \cite[\S\,3.4]{resmerita_anderssen2007}.
\begin{lemma}\label{lem:entropy_propr} 
    \begin{enumerate}
        \item \label{f_convex} 
            The function $f$ is convex.
        \item \label{f_wlsc}
            The function $f$ is weakly lower semicontinuous in $L^1(\Omega)$. 
        \item \label{f_sublevel}
            For any $c>0$, the sublevel set
            \[
                \set{v\in L^1_+(\Omega)}{f(v)\leq c}
            \] is convex, weakly closed, and weakly compact in $L^1(\Omega)$. 
        \item \label{f_dom}
            The domain of the function $f$ is strictly included in ${{L}}_+^1(\Omega)$.
        \item \label{f_intdom}
            The interior of the domain of the function $f$ is empty.
        \item \label{f_domdf}
            The set $\partial f(u)$ is nonempty if and only if $u$ belongs to $L_+^{\infty}(\Omega)$ and is bounded away from zero. In this case, $\partial f(u)=\{1+\log u\}$.
        \item \label{f_df}
            The directional derivative of the function $f$ is given by
            \[
                f'(u;v)=\int_{\Omega} v(t)[1+\log u(t)]\,dt,
            \]
            whenever it is finite.
    \end{enumerate}
\end{lemma}
Based on \cref{lem:entropy_propr}\,\ref{f_domdf}, we define in the following
\begin{equation}
    \dom\partial f=\set{u\in L^\infty_+(\Omega)}{u\text{ bounded away from zero a.e.}}.
\end{equation}

The \emph{Kullback--Leibler divergence}, which coincides with the \emph{Bregman distance} with respect to the
Bolzmann--Shannon entropy, can be defined as $d:\dom f\times \dom f\rightarrow [0,+\infty]$ by
\begin{equation}
    d(v,u)=f(v)-f(u)-f'(u;v-u),
\end{equation}
where $f'(u;\cdot)$ is the directional derivative at $u$. One can also write
\begin{equation}\label{eq:kl}
    d(v,u)=\int \left[v(t)\log\frac{v(t)}{u(t)}-v(t)+u(t) \right]\,dt,
\end{equation}
when $d(v,u)$ is finite.
We list below several properties of the Kullback--Leibler divergence.
\begin{lemma}\label{lem:kl_propr} 
    \begin{enumerate}
        \item \label{d_convex}
            The function $(v,u)\mapsto d(v,u)$ is convex.
        \item \label{d_wlsc}
            The function $d(\cdot,u^*)$ is weakly lower semicontinuous in $L^1(\Omega)$ whenever $u^*\in \dom f$.
        \item \label{d_sublevel}
            For any $c>0$ and any non-negative $u\in L^1(\Omega)$, the sublevel set \[
                \set{v\in L^1_+(\Omega)}{d(v,u)\leq c}
            \]
            is convex, weakly closed, 
            and weakly compact in $L^1(\Omega)$.
        \item \label{d_dd}
            The set $\partial d(\cdot,u^*)(u)$ is nonempty for $u^*\in \dom f$ if and only if $u$ belongs to $L_+^{\infty}(\Omega)$ and is bounded away from zero. Moreover, $\partial d(\cdot,u^*)(u)=\{\log u-\log u^*\}$.
    \end{enumerate}
\end{lemma}
Finally, the Kullback--Leibler divergence provides a bound on the $L^1$ distance.
\begin{lemma}\label{f_ineq}
    For any $u,v\in dom\,f$, one has
    \begin{equation}\label{eq:inequality}
        \|u-v\|_1^2\leq
        \left(\frac{2}{3}\|v\|_1+\frac{4}{3}\|u\|_1\right)d(v,u).
    \end{equation}
\end{lemma}

\section{Variational methods}\label{sec:variational}

Tikhonov regularization with additional convex constraints such as the non-negative cone is now classical; we refer, e.g., to \cite{Neubauer88} for linear and \cite{Chavent94,Neubauer89} for nonlinear inverse problems in Hilbert spaces and to \cite{Flemming11} for the Banach space setting. 
We therefore focus in this section on methods that are based on minimizing some combination of the regularization functional $\cR\in\{f,d(\cdot,u_0)\}$, where $u_0\in \dom f\subseteq L^1_+(\Omega)$ is an a priori guess, with the residual norm either as a penalty, i.e., as Tikhonov regularization
\begin{equation}\label{Tikhonov}
    \min_{u\in L^1_+(\Omega)} \tfrac12\|Au-\ydel\|_Y^2+\alpha \cR(u)
\end{equation}
for some regularization parameter $\alpha>0$, or as a constraint, i.e., as Morozov regularization
\begin{equation}\label{Morozov}
    \min_{u\in L^1_+(\Omega)} \cR(u)\quad\text{s.t.}\quad\|Au-\ydel\|_Y\leq\delta,
\end{equation}
where $\delta$ is the noise level according to \eqref{eq:delta}.
Throughout this section we will set $X=L^1(\Omega)$ and assume $A:X\to Y$ to be a bounded linear operator mapping into some Banach space $Y$.
Moreover, we will assume existence of a solution $u^\dag$ to \eqref{eq:op_equation} with finite entropy $\cR(u^\dag)<\infty$ (which is therefore in particular non-negative).

\subsection{Morozov-entropy regularization}

The historically first study of regularizing properties of such methods can be found in \cite{amato_hughes1991} for the Morozov-entropy method
\begin{equation}\label{eq:Morozov-entropy}
    \min_{u\in L^1_+(\Omega)} f(u)\quad\text{s.t.}\quad\|Au-\ydel\|_Y\leq \delta .
\end{equation}
The reader is referred also to \cite[Theorem 3.1]{Borwein93}, where a version of Morozov regularization is discussed.

We first of all state existence of a solution to \eqref{eq:op_equation} that maximizes the entropy, i.e., minimizes $f$. 
\begin{theorem}[existence of maximum entropy solution for exact data 
    \protect{\cite[Thm.~4]{amato_hughes1991}}]
    \label{th:maxent-exact}
    There exists a minimizer $\udag\in L^1_+(\Omega)$ of 
    \begin{equation}
        \min_{u\in L^1_+(\Omega)} f(u)\quad\text{s.t.}\quad Au=y.
    \end{equation}
\end{theorem}
\begin{theorem}[existence of regularizer \protect{\cite[Thm.~1]{amato_hughes1991}}]\label{th:Morozov-entropy-minimizer}
    For every $\delta>0$ and $\ydel\in Y$ satisfying \eqref{eq:delta}, there exists a minimizer $u^\delta$ of \eqref{eq:Morozov-entropy}.
\end{theorem}
Both theorems follow from weak compactness of sublevel sets and weak lower semicontinuity of $f$ (\cref{lem:entropy_propr}\,\ref{f_wlsc}, \ref{f_sublevel}) together with weak closedness and nonemptiness of the sets $\set{u\in L^1_+(\Omega)}{Au=y}$ and $\set{u\in L^1_+(\Omega)}{\|Au-\ydel\|_Y\leq \delta}$.

Finally, it can be shown that \eqref{eq:Morozov-entropy} indeed defines a regularization method. 
\begin{theorem}[convergence as $\delta\to0$ 
    \protect{\cite[Thm.~5]{amato_hughes1991}}]
    \label{th:Morozov-entropy-convergence}
    Let $(\ydel)_{\delta>0}$ be a family of data satisfying \eqref{eq:delta}. Then 
    \begin{equation}\label{eq:Morozov-entropy-convergence}
        \|u^\delta-\udag\|_1\to0 \text{ as }\delta\to0.
    \end{equation}
\end{theorem}

Stability of \eqref{eq:Morozov-entropy} in the sense that small perturbations in $\ydel$ lead to small perturbations in $u^\delta$ has not been shown in \cite{amato_hughes1991}. Since one is actually interested in approaching $\udag$ rather than $u^\delta$, such stability results might be considered as of lower importance than the convergence in \eqref{eq:Morozov-entropy-convergence}. We will therefore not state such stability results for the remaining variational methods either. 

Convergence rates are not stated in \cite{amato_hughes1991}, but they could be proved as in \cref{th:Tikhonov-entropy-rate} or \cref{th:Tikhonov-KL-rate} below under a similar source condition.

\subsection{Tikhonov-entropy regularization}

The work \cite{amato_hughes1991} also mentions the Tikhonov-entropy regularization
\begin{equation}\label{eq:Tikhonov-entropy0}
    \min_{u\in L^1_+(\Omega)} \tfrac12\|Au-\ydel\|_Y^2+\alpha f(u),
\end{equation}
pointing out from \cite{TikhonovArsenin77} that a regularization parameter choice $\alpha=\alpha(\delta)$ exists such that minimizers of \eqref{eq:Tikhonov-entropy0} coincide with minimizers of \eqref{eq:Morozov-entropy}. Hence, \cref{th:Morozov-entropy-convergence} also yields convergence of solutions of \eqref{eq:Tikhonov-entropy0}; however, this does not provide a concrete rule for choosing $\alpha$.

A more detailed analyis of the constrained Tikhonov-entropy regularization
\begin{equation}\label{eq:Tikhonov-entropy}
    \min_{u\in \mathcal{D}} \tfrac12\|Au-\ydel\|_Y^2+\alpha f(u)
\end{equation}
with an appropriate subset $\mathcal{D}$ of $L^1_+(\Omega)$ -- including a priori regularization parameter choice rules -- can be found in \cite{engl_landl1993}. The analysis relies on a nonlinear transformation $T:L^2_{-e^{-1/2}}(\Omega)\to L^1(\Omega)$ for $L^2_{-e^{-1/2}}(\Omega)=\set{v\in L^2(\Omega)}{v\geq {-e^{-1/2}}\text{ a.e.}}$ satisfying
\begin{equation}
    f(T(v))=\|v\|_2^2-e^{-1/2}.
\end{equation}
This leads to replacing \eqref{eq:op_equation} with the nonlinear problem $F(v)=y$
for $F:=A\circ T:L^2_{-e^{-1/2}}(\Omega)\supseteq\mathcal{B}\to Y$.
The theory on Tikhonov regularization for nonlinear problems in Hilbert spaces, in combination with a proof of weak sequential closedness of $F$, allows the authors of \cite{engl_landl1993} to prove well-definedness and convergence of minimizers of \eqref{eq:Tikhonov-entropy} under the assumption that $\mathcal{B}=T^{-1}(\mathcal{D})$ is compact in measure.
\begin{theorem}[existence of minimizers]\label{th:Tikhonov-entropy-minimizer}
    For every $\alpha>0$, $\delta>0$ and $\ydel\in Y$ satisfying \eqref{eq:delta} there exists a minimizer $u_\alpha^\delta$ of \eqref{eq:Tikhonov-entropy}.
\end{theorem}
(This result also follows from \cref{lem:entropy_propr}\,\ref{f_wlsc}, \ref{f_sublevel} together with \cite[Lem.~3.1,3.2]{engl_landl1993}.)
\begin{theorem}[convergence as $\delta\to0$ 
    \protect{\cite[Thm.~3.7]{engl_landl1993}}]\label{th:Tikhonov-entropy-convergence}
    Let $(\ydel)_{\delta>0}$ be a family of data satisfying \eqref{eq:delta}, and let $\alpha=\alpha(\delta)$ be chosen such that 
    \begin{equation}\label{eq:alpha_convergence}
        \alpha\to0 \quad\text{and}\quad\frac{\delta^2}{\alpha}\to0 \quad\text{as }\delta\to0.
    \end{equation}
    Then 
    \begin{equation}
        \|u_\alpha^\delta-\udag\|_1\to0 \quad\text{as }\delta\to0.
    \end{equation}
\end{theorem}
The reader is referred also to \cite[Theorem 3.3]{Borwein93} for convergence in the case of exact data.

Furthermore, a classical result on convergence rates for nonlinear Tikhonov regularization in Hilbert spaces from \cite{EKN89} can be employed to yield the following statement.
\begin{theorem}[convergence rates \protect{\cite[Thm.~3.8]{engl_landl1993}}]\label{th:Tikhonov-entropy-rate}
    Assume that $\udag$ satisfies the source condition
    \begin{equation}\label{eq:sourcecond_T-e0}
        1+\log \udag = A^*w
    \end{equation}
    for some sufficiently small $w\in Y^*$.
    Let $(\ydel)_{\delta>0}$ be a family of data satisfying \eqref{eq:delta}, and let $\alpha=\alpha(\delta)$ be chosen such that 
    \begin{equation}\label{eq:alpha_rate}
        \alpha\sim\delta \quad\text{as }\delta\to0.
    \end{equation}
    Then 
    \begin{equation}
        \|u_\alpha^\delta-\udag\|_1=\mathcal{O}(\sqrt{\delta}) \quad\text{as }\delta\to0.
    \end{equation}
\end{theorem}

\medskip

The work \cite{engl_landl1993} also treats the generalized negative Boltzmann--Shannon entropy
\begin{equation}\label{eq:entropy_engllandl}
    f(u)=
    \begin{cases}
        \int_{\Omega} u(t)\log \frac{u(t)}{u_0(t)}\,dt & \text{if $u\geq 0$ a.e. and $u \log \frac{u}{u_0}\in {{L}}^1(\Omega)$}, \\
        {+\infty} & \text{otherwise,}
    \end{cases}
\end{equation}
for some non-negative function $u_0$ carrying a priori information on $u$.
In this case, the source condition \eqref{eq:sourcecond_T-e0} becomes
\begin{equation}\label{eq:sourcecond_T-e}
    1+\log \frac{\udag}{u_0} = A^*w.
\end{equation}
Finally, we point out that the convergence analysis in \cite{engl_landl1993} includes the practically relevant case of inexact minimization of the Tikhonov functional.

\subsection{Tikhonov--Kullback--Leibler regularization}

Replacing the negative Boltzmann--Shannon entropy in \eqref{eq:Tikhonov-entropy} by the Kullback--Leibler divergence results in the problem
\begin{equation}\label{eq:Tikhonov-KL}
    \min_{u\in L^1_+(\Omega)} \tfrac12\|Au-\ydel\|_Y^2+\alpha d(u,u_0),
\end{equation}
which was investigated in \cite{eggermont1993}. In contrast to the analysis of the similar-looking problem \eqref{eq:Tikhonov-entropy} in \cite{engl_landl1993}, the analysis in \cite{eggermont1993} treats \eqref{eq:Tikhonov-KL} directly by convexity arguments and does not require a nonlinear transformation. 
More precisely, the fact that both parts of the Tikhonov functional can be written as Bregman distances yields the estimates
\[
    \tfrac12\|A(u-u_\alpha^\delta)\|_Y^2+\alpha d(u,u_\alpha^\delta) \leq 
    \tfrac12\|Au-\ydel\|_Y^2+\alpha d(u,u_0)
    -\tfrac12\|Au_\alpha^\delta-\ydel\|_Y^2-\alpha d(u_\alpha^\delta,u_0)
\]
for a minimizer $u_\alpha^\delta$ of \eqref{eq:Tikhonov-KL} and all $u\in L^1_+(\Omega)$; cf. \cite{Bregman67} and \cite[Thm.~3.1]{eggermont1993}. In addition, the estimate 
\[
    \tfrac12\|A(\tilde{u}_\alpha^\delta-u_\alpha^\delta)\|_Y^2+\alpha d(\tilde{u}_\alpha^\delta,u_\alpha^\delta) \leq 2 \|\tilde{y}^\delta-\ydel\|_Y^2
\]
holds for a minimizer $\tilde{u}_\alpha^\delta$ of \eqref{eq:Tikhonov-KL} with $\ydel$ replaced by $\tilde{y}^\delta$; cf. \cite[Thm.~3.3]{eggermont1993}. These two inequalities are the basis for the following results.

We first consider existence of minimizers.
Assuming the existence of a solution $u^\dag\in L^1_+(\Omega)$ to \eqref{eq:op_equation} of finite Kullback--Leibler divergence, one obtains (as in \cref{th:maxent-exact}) an existence result for exact data from the weak compactness of sublevel sets and the weak lower semicontinuity of $d(\cdot,u_0)$ (\cref{lem:kl_propr}\,\ref{d_wlsc}, \ref{d_sublevel}) as well as the weak continuity of $A$.
\begin{theorem}[existence of minimum-KL solution for exact data]
    \label{th:minKL-exact}
    For every $u_0\in\dom f$, there exists a minimizer $\udag$ of 
    \begin{equation}
        \min_{u\in L^1_+(\Omega)} d(\cdot,u_0)\quad\text{s.t.}\quad Au=y.
    \end{equation}
\end{theorem}
Next we consider the well-definedness of \eqref{eq:Tikhonov-KL} for noisy data, where one also obtains uniqueness as well as a uniform positivity property of the minimizer.
\begin{theorem}[existence, uniqueness and uniform positivity of regularizer 
    \protect{\cite[Thm.~2.4]{eggermont1993}}]\label{th:Tikhonov-KL-minimizer}
    For every $\alpha>0$, $\delta>0$, and $\ydel\in Y$ satisfying \eqref{eq:delta}, there exists a unique minimizer $u_\alpha^\delta\in L^1_+(\Omega)$ of \eqref{eq:Tikhonov-KL}.
    Moreover, $\frac{u_\alpha^\delta}{u_0}$ is bounded away from zero.
\end{theorem}
From this, the following convergence and convergence rate results can be stated.
\begin{theorem}[convergence as $\delta\to0$ 
    \protect{\cite[Thm.~4.1]{eggermont1993}}]\label{th:Tikhonov-KL-convergence}
    Let $(\ydel)_{\delta>0}$ be a family of data satisfying \eqref{eq:delta}, and let $\alpha=\alpha(\delta)$ be chosen according to \eqref{eq:alpha_convergence}. Then
    \begin{equation}\label{eq:Tikhonov-KL-convergence}
        \|u_\alpha^\delta-\udag\|_1\to0 \quad\text{as }\delta\to0.
    \end{equation}
\end{theorem}
\begin{theorem}[convergence rates; 
    \protect{\cite[Thm.~4.2]{eggermont1993}}]\label{th:Tikhonov-KL-rate}
    Assume that $\udag$ satisfies the source condition
    \begin{equation}\label{eq:sourcecond_T-KL}
        \log \frac{\udag}{u_0} = A^*w.
    \end{equation}
    Let $(\ydel)_{\delta>0}$ be a family of data satisfying \eqref{eq:delta}, and let $\alpha=\alpha(\delta)$ be chosen according to \eqref{eq:alpha_rate}.
    Then
    \begin{equation}\label{eq:rateTikhonovKL}
        \|u^\delta-\udag\|_1=\mathcal{O}(\sqrt{\delta}) \quad\text{as }\delta\to0.
    \end{equation}
\end{theorem}
We remark that in \cite{eggermont1993}, compactness of $A$ is assumed, but an inspection of the proofs shows that actually boundedness of $A$ suffices, since this implies weak lower semicontinuity of the mapping $u\mapsto \tfrac12\|Au-\ydel\|_Y^2$ and therefore (by \cref{lem:kl_propr}\,\ref{d_wlsc} and the elementary inequality $\liminf(a_n)+\liminf(b_n)\leq\liminf(a_n+b_n)$) of the Tikhonov functional.

\subsection{Nonquadratic data misfit}

In the remainder of this section, we remark on some possible extensions and open problems. 

First, the quadratic term $\frac12\|Au-\ydel\|_Y^2$ can be replaced by some other convex data misfit functional to take into account special features of the data or of the measurement noise.
In particular, we mention the case of non-negative data resulting from a positivity preserving operator $A:L^1_+(\Omega)\to L^1_+(\Omega)$ as in, e.g., \cite{resmerita_anderssen2007}) or Poisson noise as in, e.g., \cite{WernerHohage12}.

Indeed, it was shown in \cite{resmerita_anderssen2007} that \cref{th:Tikhonov-KL-minimizer,th:Tikhonov-KL-convergence,th:Tikhonov-KL-rate} also hold for the entropy--entropy regularization 
\begin{equation}\label{eq:Tikhonov-KL-KL}
    \min_{u\in L^1_+(\Omega)} d(\ydel,Au)+\alpha d(u,u_0)
\end{equation}
with $Y=L^1(\Omega)$ and \eqref{eq:delta} replaced by 
\begin{equation}\label{eq:delta_d}
    d(\ydel,y)\leq\delta^2,
\end{equation}
provided $A$ is positivity preserving in the sense that $x>0$ almost everywhere implies $Ax>0$ almost everywhere.

The key estimates for proving convergence and convergence rates in this case are
\[ 
    d(\ydel,Au_\alpha^\delta)+\alpha d(u_\alpha^\delta,u_0)
    \leq d(\ydel,A\udag)+\alpha d(\udag,u_0)
    \leq \delta^2+\alpha d(\udag,u_0),
\]
which follows from minimality and \eqref{eq:delta_d}, and
\begin{equation}
    \begin{aligned}
        d(\ydel,Au_\alpha^\delta)+\alpha d(u_\alpha^\delta,\udag)
        &= d(\ydel,Au_\alpha^\delta)\\
        \MoveEqLeft[-1]+\alpha \left(d(u_\alpha^\delta,\udag)-d(\udag,u_0)-\int_\Omega\log\frac{\udag(t)}{u_0(t)}(u_\alpha^\delta(t)-\udag(t))\, dt\right)\\
        &\leq\delta^2-\alpha\langle w, A (u_\alpha^\delta-\udag)\rangle_{Y^*,Y}\\
        &=\delta^2+\alpha\langle w, y-\ydel\rangle_{Y^*,Y}-\alpha\langle w, A u_\alpha^\delta-\ydel\rangle_{Y^*,Y}\\
        &\leq\delta^2+\tfrac43\alpha\|w\|_{Y^*}(\|\ydel\|_1+\|y\|_1+\|Au_\alpha^\delta\|_1)^{\frac12}
        (\delta+d(\ydel,Au_\alpha^\delta)^{\frac12}).
    \end{aligned}
\end{equation}
The last two inequalities hold due to the source condition \eqref{eq:sourcecond_T-KL} and to \eqref{eq:inequality}.
By using the a priori choice $\alpha\sim\delta$, the latter estimate implies $d(\ydel,Au_\alpha^\delta)=\mathcal{O}(\delta^2)$ and 
$d(u_\alpha^\delta,\udag)=\mathcal{O}(\delta)$. The inequality \eqref{eq:inequality} now yields the rate \eqref{eq:rateTikhonovKL}.

Finally, as mentioned in \cite{resmerita_anderssen2007}, convergence can also be extended to the symmetric Kullback-Leibler functional as a regularizing term in
\begin{equation}\label{eq:Tikhonov-KL-KL-symm}
    \min_{u\in L^1_+(\Omega)} d(\ydel,Au)+\alpha (d(u,u_0)+d(u_0,u)).
\end{equation}

\subsection{Measure space solutions}

In particular in the context of probability densities, it could be appropriate to look for solutions in the space of positive measures $\mathcal{M}_+(\Omega)$ instead of $L^1_+(\Omega)$, i.e., consider
\begin{equation}\label{Tikhonov_measure}
    \min_{u\in \mathcal{M}_+(\Omega)} \tfrac12\|Au-\ydel\|_Y^2+\alpha \cR(u).
\end{equation}
For a definition of entropy functionals on measure spaces we refer, e.g., to \cite{StummerVajda12}.
Minimization of some data misfit with a norm penalty but without imposing non-negativity as in \cite{BrediesPikkarainen13,ClasonKunisch12} or with non-negativity constraints but without adding a penalty as in \cite{ClasonSchiela17} has been shown to yield sparse solutions in the sense that the support of the minimizers will typically have zero Lebesgue measure.
Here it would be interesting to investigate whether an entropy penalty $\cR$ could overcome singularity of the optimality conditions (cf., e.g., \cite{ClasonSchiela17}) for attainable data that is present also for the measure space norm as penalty. Other relevant possibilities include the Wasserstein-1 and Kantorovich--Rubinstein norms considered in \cite{Lorenz14}.

\subsection{Nonlinear Problems}

A natural extension would be to consider Tikhonov regularization for a nonlinear operator $F:\dom(F)\to Y$, i.e., 
\begin{equation}
    \min_{u\in \mathcal{D}} \tfrac12\|F(u)-\ydel\|_Y^2+\alpha \cR(u),
\end{equation}
for $\cR\in\{f,d(\cdot,u_0)\}$ and $\mathcal{D}\subseteq \dom(F)\cap L^1_+(\Omega)$. For Tikhonov-entropy regularization, the analysis of \cite{engl_landl1993} by way of nonlinear transformation was extended to nonlinear operators in \cite{Landl96}.
In addition, recent analysis for Tikhonov regularization with abstract regularization functionals from, e.g., \cite{HKPS07,DissFlemming,DissPoeschl,DissWerner} together with \cref{lem:entropy_propr,lem:kl_propr} shows that the existence and convergence results from \cref{th:Tikhonov-entropy-minimizer,th:Tikhonov-entropy-convergence,th:Tikhonov-KL-minimizer,th:Tikhonov-KL-convergence} remain valid if $A$ is replaced by a nonlinear operator which is weakly continuous from $L^1(\Omega)$ to $Y$ and if $\mathcal{D}$ is weakly closed in $L^1(\Omega)$.
Furthermore, \cref{th:Tikhonov-entropy-rate,th:Tikhonov-KL-rate} can be recovered in the nonlinear case by replacing $A$ in the source conditions \eqref{eq:sourcecond_T-e0}, \eqref{eq:sourcecond_T-KL} by the Fréchet derivative $F'(\udag)$. 
However, uniform positivity results like those from \cref{th:Tikhonov-KL-minimizer} and \cite[Sec.~4.2]{resmerita_anderssen2007} do not follow from this theory and would have to be subject of additional investigations. 

\subsection{Ivanov regularization}

Ivanov regularization (also called method of quasi-solutions), cf., e.g., \cite{DombrovskajaIvanov65, Ivanov62, Ivanov63, IvanovVasinTanana02, LorenzWorliczek13, NeubauerRamlau14, SeidmanVogel89}, 
defines $u_\rho^\delta$ as a solution to 
\begin{equation}\label{eq:Ivanov}
    \min_{u\in L^1_+(\Omega)} \|Au-\ydel\|_Y \quad\text{s.t.}\quad \cR(u)\leq\rho 
\end{equation}
with $\cR\in\{f,d(\cdot,u_0)\}$. Here, regularization is controlled by the parameter $\rho>0$, with larger parameters corresponding to weaker regularization. If the respective minimizers are unique,
all three variational regularization methods (Tikhonov, Morozov, and Ivanov) are equivalent for a certain choice of the regularization parameters $\alpha$ and $\rho$, cf. \cite[Thm.~2.3]{LorenzWorliczek13}. Nevertheless, a practically relevant regularization parameter choice might lead to different solutions; the three formulations also entail different numerical approaches, some of which might be better suited than others in concrete applications. Furthermore, as a counterexample in \cite{LorenzWorliczek13} shows, the methods are no longer equivalent in the non-convex case arising from a nonlinear forward operator $F$. 

Concerning well-definedness and convergence for \eqref{eq:Ivanov}, one can again rely on general results on the entropy functionals as stated in \cref{sec:preliminaries}.
Indeed, the properties of sublevel sets according to \cref{lem:entropy_propr}\,\ref{f_sublevel} or \cref{lem:kl_propr}\,\ref{d_sublevel} together with weak sequential lower semicontinuity of the mapping $u\mapsto\frac12\|Au-\ydel\|^2$ guarantee existence of a minimizer for any $\rho>0$.
In case the maximal entropy or the minimal Kullback--Leibler divergence of a solution to \eqref{eq:op_equation} is known, the ideal parameter choice is of course $\rho=\cR(\udag)$; note that this choice of $\rho$ is independent of $\delta$. 
In this case, we obtain from minimality of $u_\rho^\delta$ that
\[
    \|Au_\rho^\delta-\ydel\|_Y\leq \|A\udag-\ydel\|_Y\leq\delta.
\]
We can thus argue similarly to the convergence proof for the Morozov formulation \eqref{eq:Morozov-entropy} to obtain convergence $\|u_\rho^\delta-\udag\|_1\to0$ as $\delta\to0$. Convergence rates under source conditions of the type \eqref{eq:sourcecond_T-e0} can also be derived. To see this in case $\cR=f$, observe that minimality of $u_\rho^\delta$ and admissibility of $\udag$ together with \cref{lem:entropy_propr}\,\ref{f_df} and \eqref{eq:sourcecond_T-e0} yields 
\[
    d(u_\rho^\delta,\udag)=f(u_\rho^\delta)-f(\udag)-f'(\udag,u_\rho^\delta-\udag)
    \leq 
    -\langle w,A(u_\rho^\delta-\udag)\rangle_{Y^*,Y}\leq 2\|w\|_{Y^*}\delta.
\]
Hence, \eqref{eq:inequality} implies the rate \eqref{eq:rateTikhonovKL}.
A practical choice of $\rho$ can be carried out, e.g., by Morozov's discrepancy principle, see \cite{KK17}. 

\section{Iterative methods}\label{sec:iterative}

In practice, solutions to the approaches given in \cref{sec:variational} cannot be computed directly but require iterative methods. This makes applying iterative regularization methods directly to \eqref{eq:op_equation} attractive. For the particular case of non-negativity constraints, there are two general approaches: The constraints can be imposed during the iteration by projection, or the iteration can be constructed such that non-negativity of the starting value is preserved. In this section, we will review examples of both methods.

\subsection{Projected Landweber method for non-negative solutions of linear ill-posed equations}

The classical Landweber method for the solution of \eqref{eq:op_equation} in Hilbert spaces consists in choosing $u_0=0$, $\tau \in(0, 2\norm{A}^{-2})$, and setting
\begin{equation}
    u_{k+1} = u_k + \tau A^*(y-Au_k),\qquad k=0,\dots.
\end{equation}
By spectral methods, one can show that the iterates converge strongly to the minimum norm solution $u^\dag$ for exact data $y\in \ran A$.
For noisy data $y=y^\delta \in\overline{\ran A}\setminus\ran A$, one also initially observes convergence, but at some point the iterates start to diverge from the solution; this behavior is often referred to as \emph{semi-convergence}. It is therefore necessary to choose an appropriate stopping index $k_*:=k_*(\delta,y^\delta)<\infty$ such that $u_{k_*}^\delta \to u^\dag$ as $\delta \to 0$; a frequent choice is a discrepancy principle, e.g., of Morozov.

This method was generalized in \cite{Eicke92} to constrained inverse problems of the form
\begin{equation}
    Au = y \quad\text{s.t.}\quad u\in C
\end{equation}
for a convex and closed set $C\subset X$; in our context, the obvious choice is $X=L^2(\Omega)$ and 
\begin{equation}
    C = \left\{u\in X: u(x)\geq 0 \text{ for almost every }x\in \Omega\right\}.
\end{equation}
The corresponding \emph{projected Landweber method} then consists in the iteration
\begin{equation}\label{eq:projland}
    u_{k+1} = P_C\left[u_k + \tau A^*(y-Au_k)\right],\qquad k=0,\dots,
\end{equation}
where $P_C$ denotes the metric (in our case pointwise almost everywhere) projection onto $C$.
This coincides with a \emph{forward--backward splitting} or \emph{proximal gradient descent} applied to $\norm{Au-y}_X^2 + \delta_C(u)$, where $\delta_C$ denotes the indicator function of $C$ in the sense of convex analysis; see, e.g., \cite{Combettes05}. Thus, a standard proof yields weak convergence of the iterates in the case of exact data $y\in A(C):=\set{Au}{u\in C}$.
\begin{theorem}[\protect{\cite[Thm.~3.2]{Eicke92}}]
    Let $u_0=0$ and $\tau \in(0, 2\norm{A}^{-2})$. If $y\in A(C)$, then the sequence of iterates $\{u_k\}_{k\in\N}\subset C$ of \eqref{eq:projland} converges weakly to a solution $u^\dag\in C$ of $Au=y$.
\end{theorem}
In contrast to the unconstrained Landweber iteration, strong convergence can only be shown under additional restrictive conditions or by including additional terms in the iteration; in the setting considered here, this holds if $\Id - \tau A^*A$ is compact, see \cite[Thm.~3.3]{Eicke92}.

Regarding noisy data $y^\delta\notin\ran A$, the following stability estimate holds.
\begin{theorem}[\protect{\cite[Thm.~3.4]{Eicke92}}]
    Let $u_0^\delta=u_0=0$ and $\tau \in(0, 2\norm{A}^{-2})$. If $y\in A(C)$ and $y^\delta \in Y$ with $\norm{y^\delta - y}_Y\leq \delta$, then the sequences of iterates $\{u_k\}_{k\in\N}\subset C$ and $\{u_k^\delta\}_{k\in\N}\subset C$ of \eqref{eq:projland} with $y$ and $y^\delta$, respectively, satisfy
    \begin{equation}\label{eq:projland_stab}
        \norm{u_k^\delta-u_k}_X \leq \tau\norm{A} \delta k,\qquad k=0,\dots.
    \end{equation}
\end{theorem}
By usual arguments, this can be used -- together with a monotonicity property for noisy data -- to derive stopping rules and thus regularization properties. However, to the best of our knowledge, this has not been done in the literature so far. It should also be pointed out that the estimate \eqref{eq:projland_stab} is weaker than in the unconstrained case, where an $\mathcal{O}(\sqrt{k})$ estimate can be shown. 

In addition, \cite{Eicke92} proposes a ``dual'' projected Landweber iteration: Setting $w_0=0$, compute for $k=0,\dots,$ the iterates
\begin{equation}
    \left\{
        \begin{aligned}
            u_k &= P_C A^*w_k,\\
            w_{k+1} &= w_k + \tau (y- Au_k).
        \end{aligned}
    \right.
\end{equation}
(This can be interpreted as a \emph{backward--forward splitting}.)
Under the same assumptions as above, one obtains strong convergence $u_k\to u^\dag$ (without assuming compactness of $\Id - \tau A^*A$), see \cite[Thm.~3.5]{Eicke92}, and the stability estimate \eqref{eq:projland_stab}, see \cite[Thm.~3.6]{Eicke92}. Numerical examples for integral equations with non-negative solutions in $L^2(\Omega)$ using both methods can be found in \cite{Eicke92}. Acceleration by stationary preconditioning -- i.e., replacing the scalar $\tau$ by a fixed self-adjoint, positive definite, linear operator $D$ -- was considered in \cite{Piana97}. It is an open problem whether further acceleration by Nesterov-type extrapolation or a more general inertial approach is possible.

\bigskip

If $X$ and $Y$ are Banach spaces, the above iterations are not applicable. A version of Landweber iteration in Banach spaces that can treat convex constraints has been proposed in \cite{BotHein12}. The iteration can be formulated as
\begin{equation}
    \left\{
        \begin{aligned}
            \xi_{k+1} &= \xi_k - \tau A^* J_Y (y-Ax_k),\\
            x_{k+1} &\in \arg\min_{x\in X} G(x) - \dual{\xi_{k+1},x}_X,
        \end{aligned}
    \right.
\end{equation}
for $x_0=0$ and $\xi_0=0$, where $J_Y$ denotes the so-called duality mapping between $Y^*$ and $Y$, and $G:X\to\R\cup\{\infty\}$ is proper, convex and lower semicontinuous. Here the choice $G=\delta_C$ yields the projected Landweber iteration in Banach spaces, while $G=\delta_C + f$ would be the choice if a non-negative minimum-entropy solution is searched for. However, convergence in \cite{BotHein12} could only be shown under the assumption that the interior of $C$ is non-empty in $X$ (which is not the case for $X=L^p(\Omega)$, $p<\infty$; see, e.g., \cite{BorweinLimber96}), $G$ is $p$-convex with $p\geq 2$ (in particular, excluding both $G=\delta_C$ and $G=\delta_C+f$), and $Y$ is uniformly smooth (requiring $Y$ to be reflexive). The first assumption was removed in \cite{Jin13}, allowing application to the case $X=L^p(\Omega)$ for $2\leq p<\infty$ with $G=\delta_C + f + \frac1p\norm{\cdot}_X^p$ and $Y=L^q(\Omega)$ for $2\leq q<\infty$. 
Another open issue is the practical realization for $p,q>2$, in particular of the second step, which requires computing a generalized metric projection in a Banach space.

Alternatively, a natural first step of moving from Tikhonov-type variational regularization towards iterative methods is the so-called non-stationary Tikhonov regularization \cite{hanke_groetsch1998} (which is a proximal point method), whose entropy-based version can be formulated as 
\begin{equation}\label{eq:Tikhonov-nonst-entropy}
    u_{k}=\text{arg}\min_{u\in L^1_+(\Omega)} \tfrac12\|Au-\ydel\|_Y^2+\alpha_k d(u,u_{k-1}),
\end{equation}
where $\{\alpha_k\}_{k\in\N}$ is a bounded sequence of positive numbers. This has been shown to converge in finite dimension; see, e.g., \cite{iusem1995} and the references therein. We expect that an analysis of the infinite-dimensional counterpart can be carried out using the tools presented in this review.

\subsection{EM method for integral equations with non-negative data and kernel}\label{sect:em}

We now consider \eqref{eq:op_equation} in the special case that $A$ is a Fredholm integral operator of the first kind, i.e.,
\begin{equation}\label{eq:integral_operator}
    A:L^1(\Omega)\rightarrow L^1(\Sigma),\qquad   (Au)(s)=\int_{{\Omega}} a(s,t)u(t)\,dt\,,
\end{equation}
where ${\Omega},{\Sigma}\subset\R^d$, $d\geq 1$, are compact, and the kernel $a$ and the data $y$ are positive pointwise almost everywhere. In this case, the following multiplicative iteration can be seen to preserve non-negativity for $u_0\geq0$:
\begin{equation}\label{eq:em}
    u_{k+1}(t)=u_k(t){\int}_{\Sigma}\frac{a(s,t)y(s)}{(Au_k)(s)}\,ds,\quad t\in\Omega,\qquad k=0,\dots.
\end{equation}
This method was introduced in \cite{kondor1983} as the method of convergent weights, motivated by some problems arising in nuclear physics.
Writing this concisely as
\begin{equation}\label{eq:em1}
    u_{k+1}=u_k\,A^*\frac{y}{Au_k},\qquad k=0,\dots,
\end{equation}
where the multiplication and division are to be understood pointwise almost everywhere, relates \eqref{eq:em} to the popular method known in the finite-dimensional setting as the expectation-maximization (EM) algorithm for Poisson models for PET, cf. \cite{shepp_vardi82, vardi_shepp_kaufmann85}, and as the Lucy--Richardson algorithm in astronomical imaging, see \cite{richardson72, lucy75}.

The study of \eqref{eq:em} was initiated by the series of papers \cite{muelthei_schorr1987, muelthei_schorr1989, muelthei1992} primarily for the setting $A:C([0,1])\to C([0,1])$. More precisely, some monotonicity features have been derived, while convergence has not been shown yet.
Modified EM algorithms allowing for better convergence properties have been investigated in infinite-dimensional setting; see \cite{eggermont_lariccia1996,eggermont1999}.

In the following, we summarize, based on \cite{iusem1992,eggermont1999}, the convergence properties for the case $A:L^1(\Omega)\to L^1(\Sigma)$ using a similar notation as in \cite{resmerita_engl_iusem2007, resmerita_engl_iusem2008}.
Specifically, we make the following assumptions:
\begin{enumerate}[label=(A\arabic*),ref=A\arabic*]
    \item\label{ass:em1}
        The kernel $a$ is a positive and measurable function satisfying
        \begin{equation}
            \int_{\Sigma} a(s,t)\,ds=1 \qquad\text{for almost all } t\in \Omega.
        \end{equation}

    \item\label{ass:em2}
        There exist $m, M>0$ such that
        \begin{equation}
            m\leq a(s,t)\leq M \qquad\text{a.e. on }\Sigma\times\Omega.
        \end{equation}

    \item\label{ass:em3}
        The exact data $y$ in \eqref{eq:op_equation} satisfies $\int_{\Sigma}y(s)\,ds=1$ and
        \begin{equation}
            y(s)\leq M' \qquad\text{a.e. on }\Sigma.
        \end{equation}
        for some $M'>0$.

    \item\label{ass:em4}
        Equation \eqref{eq:op_equation} admits a solution $u^\dag\in L^1_+(\Omega)\setminus\{0\}$.
\end{enumerate}
Furthermore, let 
\begin{equation}\label{delta}
    \Delta=\set{u\in L_+^1(\Omega)}{\int_{{\Omega}} u(t)\,dt=1\,}.
\end{equation}

By noticing that a positive solution of \eqref{eq:op_equation} is also a minimizer of the function $u\mapsto d(y,Au)$ subject to $u\geq 0$ (which is related to a maximum likelihood problem in statistical setting), we formulate a classical monotonicity result for \eqref{eq:em} for exact data.
\begin{proposition}[\protect{\cite[Prop.~3.3]{resmerita_engl_iusem2007}}]\label{p2} 
    Let \eqref{ass:em1} and \eqref{ass:em3} be satisfied and let $u_0\in\Delta$ such that $d(u^\dag,u_0)<\infty$. Then, for any $k\geq 0$, the iterates $u_k$ generated by \eqref{eq:em} satisfy
    \begin{align}
        d(u^\dag,u_k)&<\infty,\\
        d(u_{k+1}, u_k)&\leq d(y,Au_k)-d(y,Au_{k+1}),\\
        d(y,Au_k)-d(y,Au^\dag)&\leq d(u^\dag,u_k)-d(u^\dag,u_{k+1}).
    \end{align}
    Therefore, the sequences $\{d(u^\dag,u_k)\}_{k\in\mathbb{N}}$ and $\{d(y,Au_k)\}_{k\in\mathbb{N}}$ are nonincreasing. Moreover, 
    \begin{align}
        \lim_{k\rightarrow\infty}d(y,Au_k)&=d(y,Au^\dag),\\
        \lim_{k\rightarrow\infty}d(u_{k+1},u_k)&=0. 
    \end{align}
\end{proposition}
We point out that from this result, one also obtains that $\lim_{k\rightarrow\infty}\|Au_k-y\|_1=0$ and $\lim_{k\rightarrow\infty}\|u_{k+1}-u_k\|_1=0$.

Since ill-posedness is an infinite-dimensional phenomenon, and the EM algorithm has been shown to be highly unstable, we recall also the approach in \cite{resmerita_engl_iusem2007} which investigates the noise influence on the iterations. Similar properties of the iterates (as in the noise free data case) are derived there by stopping the procedure according to a discrepancy rule, as one can see below.

For the noisy data $y^\delta$, we make the following assumptions.
\begin{enumerate}[label=(A\arabic*),ref=A\arabic*,start=5]
    \item\label{ass:em5}
        $y^\delta\in L^{\infty}(\Sigma)$ satisfies $\int_{\Sigma}y^\delta(s)\,ds=1$ and
        \begin{equation}\label{eq:noise_norm}
            \|y^\delta-y\|_1\leq\delta, \qquad \delta>0;
        \end{equation}
    \item\label{ass:em6} There exist $m_1,M_1>0$ such that for all $\delta>0$,
        \begin{equation}\label{eq:away}
            m_1\leq\ydel(s)\leq M_1, \qquad\text{a.e. on }\Sigma,
        \end{equation}
\end{enumerate}
For further use, we define
\begin{equation}\label{eq:def_gamma}
    \gamma := \max \left\{ \big|{} \ln \tfrac{m_1}{M} \big| \, , \,
    \big|{}\ln \tfrac{M_1}{m} \big| \right\},
\end{equation}
where $m,M>0$ are the constants from \eqref{ass:em2}.

Let now $u_k^\delta$ denote the iterates generated by \eqref{eq:em} with $y^\delta$ in place of $y$ and $u_0^\delta = u_0\in \Delta$.
In this case, the iterates get closer and closer to the solution as long as the residual lies above the noise level.
\begin{theorem}[\protect{\cite[Thm.~6.3]{resmerita_engl_iusem2007}}]\label{t1} 
    Fix $\delta>0$. If assumptions \eqref{ass:em1}--\eqref{ass:em6} are satisfied, then 
    \begin{equation}\label{eq:mon_kl}
        d(u^\dag,\ukdn)\leq d(u^\dag, \ukd)
    \end{equation}
    for all $k\geq0$ such that
    \begin{equation}\label{rule}
        d(\ydel,A\ukd)\geq {\delta} \gamma.
    \end{equation}
\end{theorem}

The above result indicates a possible choice of the stopping index for the algorithm \eqref{eq:em} as
\begin{equation}\label{eq:stop_index}
    k_*(\delta)=\min\left\{k\in\mathbb{N}: d(\ydel,A\ukd)\leq {\tau\delta}\gamma\right\}
\end{equation}
for some fixed $\tau>1$ and $\gamma$ as given by \eqref{eq:def_gamma}.
The next statement guarantees existence of such a stopping index.
\begin{theorem}[\protect{\cite[Thm.~6.4]{resmerita_engl_iusem2007}}]\label{t2} Let assumptions \eqref{ass:em1}--\eqref{ass:em6} be satisfied and choose $u_0\in\Delta$ such that $d(u^\dag,u_0)<\infty$. Then:
    \begin{enumerate}
        \item For all $\delta>0$, there exists a $k_*(\delta)$ satisfying \eqref{eq:stop_index} and
            \begin{equation}\label{eq:finite_index}
                {k_*(\delta)\tau\delta} \gamma\leq k_*(\delta)d(\ydel,Au^\delta_{k_*(\delta)-1})\leq d(u^\dag,u_0)+{k_*(\delta)\delta} \gamma.
            \end{equation}
        \item The stopping index $k_*(\delta)$ is finite with $k_*(\delta)=O\left({\delta}^{-1}\right)$ and
            \begin{equation}\label{eq:images}
                \lim_{\delta\rightarrow 0^+}\|Au^{\delta}_{k_*(\delta)}-y\|_p=0,
            \end{equation}
            for any $p\in[1,+\infty)$.
    \end{enumerate}
\end{theorem}

An interesting open problem would be to investigate the behavior of \eqref{eq:em} in conjunction with other stopping rules, e.g., a monotone error-type rule defined by means of the KL divergence in a way similar to the one dealt with in \cite{hamarik_kaltenbacher_kangro_resmerita2016}.

\subsection{Modified EM algorithms}

In this section, we present some modifications of algorithm \eqref{eq:em} which improve its stability or its performance.

\subsubsection{EM algorithms with smoothing steps}

In order to stabilize \eqref{eq:em}, the work \cite{silverman_etal1990} proposed the so-called EMS algorithm
\begin{equation}
    u_{k+1}=S\left(u_k\,A^*\frac{y}{Au_k}\right),\qquad k=0,\dots,
\end{equation}
with $u_0\equiv 1$ and the smoothing operator
\begin{equation}
    Su(s)=\int_\Omega b(s,t)u(t)\,dt,
\end{equation}
where $b:\Omega\times\Omega\to\R$ is continuous, positive and obeys a normalization condition similar to \eqref{ass:em1}.
Note that \cite{silverman_etal1990} presents the EMS method from a statistical perspective, while the continuous formulation mentioned above can be found in \cite{eggermont1999}.
Although this yields faster convergence in practice, more information on limit points (or the unique limit point, as the numerical experiments strongly suggest) has not been provided. 

The work \cite{eggermont1999} also proposes the following nonlinear smoothing procedure, called NEMS:
\begin{equation}
    u_{k+1}=S\left(\mathcal{N}(u_k)\,A^*\frac{y}{Au_k}\right),\qquad k=0,\dots,
\end{equation}
with $u_0\equiv 1$ and 
\begin{equation}
    \mathcal{N}u(t)=\exp\big([S^*(\log u)](t)\big) \qquad\text{for all } t\in\Omega.
\end{equation}
Properties typical to EM iterations are shown in \cite{eggermont1999}. In particular, \cite[Thm.~4.1]{eggermont1999} proves that the iterates produced by the NEMS algorithm converge to solutions of
\begin{equation}
    \min_{u\in L^1_+(\Omega)} d(y,A\mathcal{N}u)-\mathcal{N}u+\int_\Omega u,
\end{equation}
in analogy with \eqref{eq:em} which is designed for approximating minimizers in $L^1_+(\Omega)$ of $d(y,Au)$.

\subsubsection{EM-Kaczmarz type algorithms}

The ordered-subsets expectation maximization algorithm (OS-EM) proposed in \cite{hudson_larkin1994} is a variation of the EM iteration that has proved to be quite efficient in computed tomography. It resembles a Kacmarz-type method, being conceptually based on grouping the data $y$ into an ordered sequence of $N$ subsets $y_j$. A single outer iteration step then is composed of $N$ EM-steps, where in each step $j$ one updates the current estimate by working only with the corresponding data subset $y_j$. An extension of the OS-EM to the infinite-dimensional setting was introduced in \cite{haltmeier_leitao_resmerita2009}, and numerical experiments reported there indicate this to be at least as efficient as the classical discrete OS-EM algorithm.

Let $\Sigma_j$ be (not necessarily disjoint) subsets of $\Sigma$ with
$ \Sigma_0 \cup \dots \cup \Sigma_{N-1}=\Sigma$, and denote $y_j := y|_{\Sigma_j}$.
Set $a_j := a|_{\Omega \times \Sigma_j}$.

Thus, one can rewrite 
\eqref{eq:integral_operator} as
\begin{equation} \label{eq:sys-ops}
    A_j: L^1(\Omega) \to L^1(\Sigma_j),\qquad (A_j u)(s)  := \int_\Omega a_j(s,t) \, u(t) \, dt,  \qquad j = 0, \dots, N-1 .
\end{equation}
Then \eqref{eq:integral_operator} can be formulated as 
a system of integral equations of
the first kind
\begin{equation} \label{eq:sys-ip}
    A_j u  =  y_j , \qquad j = 0, \dots, N-1.
\end{equation}
Clearly, $u$ is a solution of
\eqref{eq:sys-ip} if and only if $u$ solves \eqref{eq:integral_operator}.
Without loss of generality, one can work with a common domain $\Sigma$ instead of
$\Sigma_j$, thus considering $A_j: L^1(\Omega) \to L^1(\Sigma)$
and $y_j \in L^1(\Sigma)$. Therefore, the system
\eqref{eq:sys-ip} can be solved by simultaneously minimizing
\begin{equation}
    d(y_j, A_j u), \qquad j = 0, \dots, N-1 .
\end{equation}

The corresponding OS-EM algorithm for solving
system \eqref{eq:sys-ip} can be written in the form
\begin{equation} \label{eq:def-osem}
    u_{k+1}
    =
    u_k \int_\Sigma
    \frac{a_j(s,\cdot) \, y_j(s)}{(A_j u_k)(s)} \, ds,
    \qquad k=0,\dots,
\end{equation}
where $j = [k] := (k$ mod $N)$. 

Under assumptions similar to the ones in \cref{sect:em}, the following results hold for exact data.
\begin{theorem}[\protect{\cite[Thm.~3.3]{haltmeier_leitao_resmerita2009}}]
    Let the sequence $\{u_k\}_{k\in\N}$ be
    defined by iteration \eqref{eq:def-osem}, and let $u^\dag \in \Delta\setminus\{0\}$ be a  solution of \eqref{eq:sys-ops} with $d(u^\dag, u_0) < \infty$.
    Then one has
    \begin{itemize}
        \item[(i)] $f_{[k]}(u_{k+1}) \le f_{[k]}(u_k)$, for every $k =0,\dots$;
        \item[(ii)] the sequence $\{ d(u^\dag, u_k) \}_{k\in\N}$ is nonincreasing;
        \item[(iii)] $\lim\limits_{k\to\infty} f_{[k]}(u_k) = 0$;
        \item[(iv)] $\lim\limits_{k\to\infty} d(u_{k+1}, u_k) = 0$;
        \item[(v)] for each $0 \le j \le N-1$ and $p\in[1,\infty)$ we have
            \begin{equation} \label{eq:l1-conv}
                \lim_{m\to\infty} \| A_j u_{j+mN} - y_j \|_{p} = 0 ;
            \end{equation}
        \item[(vi)] if $u_0\in L^\infty(\Omega)$ and $\{u_k\}_{k\in\N}$ is bounded in 
            $L^p(\Omega)$ for some $p\in (1, \infty)$, then there exists a subsequence
            converging weakly in $L^p(\Omega)$ to a solution of \eqref{eq:sys-ip}.
    \end{itemize}
\end{theorem}

Finally, we address the case of noisy data, where the noise is allowed to be different for each group of data. Specifically, we consider $y_j^\delta \in L^1(\Sigma)$ with
\begin{equation} \label{eq:noisy-data}
    \|y_j-y_j^{\delta}\|_{1} \le \delta_j , \qquad j = 0, \dots, N-1, 
\end{equation}
and set $\delta := (\delta_0, ..., \delta_{N-1})$.
Correspondingly, the stopping rules are independently defined for each group of data.
This leads to the following loping OS-EM iteration for
\eqref{eq:sys-ip} with noisy data:
\begin{equation} \label{eq:def-osem-noise1}
    u_{k+1}^\delta =
    \begin{cases}
        u_k^\delta\int_\Sigma \frac{a_{[k]}(s,\cdot) \, y_{[k]}^\delta(\cdot)}
        {(A_{[k]} u_k^\delta)(s)} \, ds 
        & d(y_j^\delta, A_j u_k^\delta) > \tau \gamma \delta_{[k]}, \\
        u_k^\delta & \text{else},
    \end{cases} 
\end{equation}
for $\tau>1$ and $\gamma$ as defined in \eqref{eq:def_gamma}.

Under similar assumptions on kernels and data as in \cref{sect:em}, one can show analogously to \cref{t2} that there exists a finite stopping index $k_*(\delta)$, after which the iteration for all groups is terminated, and that (under additional assumptions) $A_ju_{k^*(\delta)}^\delta\to y_j$ for all $j$ and $u_{k^*(\delta)}^\delta \wkto u^\dag$ in $L^p$ as $\delta \to 0$, see \cite[Thm.~4.4]{haltmeier_leitao_resmerita2009}.

\printbibliography

\end{document}